\documentstyle[12pt]{article}
\newfont{\twlvmsb}{msbm10 scaled\magstep1}
\newfont{\ninemsb}{msbm9}
\newfont{\sixmsb}{msbm6}
\newfam\msbfam
\textfont\msbfam=\twlvmsb
\scriptfont\msbfam=\ninemsb
\scriptscriptfont\msbfam=\sixmsb
\catcode`\@=11
\def\Bbb{\ifmmode\let\next\Bbb@\else
  \def\next{\errmessage{Use \string\Bbb\space only in math mode}}\fi\next}
\def\Bbb@#1{{\Bbb@@{#1}}}
\def\Bbb@@#1{\fam\msbfam#1}
\newfont{\largeeufm}{eufm10 scaled\magstep4}
\newfont{\twlveufm}{eufm10 scaled\magstep1}
\newfont{\elveufm}{eufm10 at 11pt}
\newfont{\teneufm}{eufm10}
\newfont{\nineeufm}{eufm9}
\newfam\eufam
\textfont\eufam=\twlveufm
\scriptfont\eufam=\nineeufm
\def\frak{\ifmmode\let\next\frak@\else
\def\next{\errmessage{Use \string\frak\space only in math mode}}\fi\next}
\def\frak@#1{{\fam\eufam{{#1}}}}
\catcode`@=12

\newcommand{\be}{\begin{equation}}
\newcommand{\ee}{\end{equation}}

\newcommand{\bea}{\begin{eqnarray}}
\newcommand{\eea}{\end{eqnarray}}

\newcommand{\hst}[1]{\rule{#1}{0mm}}
 
\newcommand{\what}[1]{\widehat{\rule{#1}{0pt}}}

\newcommand{\ty}{\hspace{0.04em}}

\newcommand{\sect}[1]{\setcounter{equation}{0}\section{#1}}

\newcommand{\reff}[1]{(\ref{#1})}

\newcommand{\dis}{\displaystyle}

\newcommand{\Lam}{\Lambda}

\newcommand{\de}{\delta}

\newcommand{\ve}{\varepsilon}
\newcommand{\th}{\theta}

\newcommand{\vr}{\varrho}
\newcommand{\si}{\sigma}

\newcommand{\rd}{\partial}
\newcommand{\la}{\langle}
\newcommand{\ra}{\rangle}
\newcommand{\rar}{\rightarrow}
\newcommand{\lra}{\longrightarrow}
\newcommand{\ti}{\times}
\newcommand{\op}{\oplus}
\newcommand{\bop}{\bigoplus}
\newcommand{\ot}{\otimes}

\newcommand{\KK}{\Bbb K}
\newcommand{\ZZ}{\Bbb Z}

\newcommand{\Ze}{\{0\}}
\newcommand{\cd}{\cdot}

\newcommand{\comp}{\!\stackrel{\textstyle{\hst{0ex} \atop \circ}}{\hst{0ex}}\!}

\newcommand{\wepf}{{\:\wedge\!\!\!{}_{{}_{{}_{\scriptstyle\varepsilon}}}\:}}

\begin{document}

\begin{titlepage}

\newlength{\ppn}
\newcommand{\defboxn}[1]{\settowidth{\ppn}{#1}}
\defboxn{BONN--TH--98--nn}

\hspace*{\fill} \parbox{\ppn}{BONN--TH--98--09 \\
                              May 1998        }

\vspace{18mm}

\begin{center}
{\LARGE\bf The second cohomology of $sl(m|1)$ \\[0.5ex]
           with coefficients in its enveloping algebra \\[1.0ex]
           is trivial} \\
\vspace{11mm}
{\large M. Scheunert} \\
\vspace{1mm}
Physikalisches Institut der Universit\"{a}t Bonn \\
Nu{\ss}allee 12, D--53115 Bonn, Germany \\
\vspace{5mm}
{\large R.B. Zhang} \\
\vspace{1mm}
Department of Pure Mathematics \\
University of Adelaide \\ 
Adelaide, Australia \\
\end{center}

\vspace{15mm}
\begin{abstract}
\noindent
Using techniques developed in a recent article by the authors, it is proved
that the 2--cohomology of the Lie superalgebra $sl(m|1)$\,; $m \geq 2$\,,
with coefficients in its enveloping algebra is trivial. The obstacles in
solving the analogous problem for $sl(3|2)$ are also discussed.
\end{abstract}

\vspace{\fill}
\noindent
math.QA/9805034

\end{titlepage}

\setcounter{page}{2}
\renewcommand{\baselinestretch}{1}     %
\small\normalsize

\sect{Introduction}
%
%
The present work is a direct sequel to a recent article by the authors \cite{SZ}
dealing with the cohomology of Lie superalgebras (for a list of pertinent
references, see Ref.~\cite{SZ}). One of the main goals of these papers is to
prove or disprove, for as many basic classical simple Lie superalgebras $L$
as possible, that
 \be H^2(L,U(L)) = \Ze \,, \label{claim} \ee
where $U(L)$ is the enveloping algebra of $L$\,, endowed with the adjoint
action of $L$\,. As explained in Ref.~\cite{SZ}, Eq.~\reff{claim} implies
that the associative superalgebra $U(L)$ does not admit of any non--trivial
formal deformations in the sense of Gerstenhaber \cite{Ge}. In Ref.~\cite{SZ},
we have shown that Eq.~\reff{claim} holds for the $osp(1|2n)$ algebras and for
$sl(2|1)$. Here, we are going to prove that it is true for $sl(m|1)$ with
$m \geq 2$\,.

The setup of our paper is the following. In Sec.~2 we introduce our notation
and specialize some of the results of Ref.~\cite{SZ} to the case of present
interest, the Lie superalgebras $sl(m|1)$ with $m\geq 2$\,. Sec.~3 contains
the proof of Eq.~\reff{claim} for these algebras. A short discussion follows
in Sec.~4.

The paper is closed by an appendix, in which we consider the analogous problem
for the algebra $sl(3|2)$. Unfortunately, up to now we have not been able to
prove or disprove Eq.~\reff{claim} in this case. However, we think it may be
worthwhile to show the obstacles which one has to overcome if one wants to
proceed along the lines described in Ref.~\cite{SZ}.

We close this introduction by specifying some of our conventions. Throughout
the present work the base field $\KK$ is assumed to be {\em algebraically
closed} and to have characteristic zero. The multiplication in a Lie
superalgebra will be denoted by a pointed bracket $\la\;,\;\ra$. All modules
over a Lie superalgebra are assumed to be graded.
%

\sect{Reminder of some previous results}
%
%
Let us first explain our notation and conventions. For later use (in the
appendix) we describe them for an arbitrary special linear Lie superalgebra
$sl(m|n)$, with $m,n \geq 1$ and $m \neq n$\,. Quite generally, we follow
Ref.~\cite{SZ}, in particular, see Example 4 in Sec.~3. Thus we use the
generators
 \be  X_{ij} = E_{ij} - \frac{1}{m-n}\,\si_i\,\de_{ij}\,I \,, \label{gen} \ee
(with $i,j \in \{1,2\ty,\ldots,m+n\}$) of $sl(m|n)$, where the $E_{ij}$ are
the standard basic \\
$(m+n) \ti (m+n)$ matrices, where $I$\/ is the
$(m+n) \ti (m+n)$ unit matrix, and where
 \be  \si_i \,=\, \Bigg\{ \begin{array}{r@{\;\;\;\mbox{if}\;\;\;}l}
                                    1 & 1 \le i \le m \\
                                   -1 & m+1 \le i \le m+n \,.
                          \end{array}                               \ee
Consequently, we have  
 \be  \sum_{i=1}^{m+n}X_{ii} = 0 \,, \label{slgen}\ee
and the Cartan subalgebra $h$ of $sl(m|n)$, consisting of the diagonal matrices
in $sl(m|n)$, is spanned by the elements $X_{ii}$\,. The usual basic linear
forms on the Car\-tan subalgebra of $gl(m|n)$, associating to a diagonal matrix
its $i$\,th diagonal element, yield by restriction to $h$ the linear forms
$\ve_i \in h^{\ast}$ given by
 \be  \ve_i(X_{jj}) = \de_{ij} - \frac{1}{m-n}\,\si_j \,, \label{defeps}\ee
which satisfy the equation
 \be  \sum_{i=1}^{m+n}\si_i \,\ve_i = 0 \,. \label{releps} \ee
In terms of the $\ve_i$\,, every linear form $\Lam$ on $h$ (in particular, any
weight of an $sl(m|n)$--module) can be uniquely written in the form
 \be  \Lam = \sum_{i=1}^{m+n}L_i \,\ve_i \,, \label{defweight} \ee
with scalars $L_i$ such that
 \be  \sum_{i=1}^{m+n}L_i = 0 \,. \label{relweight} \ee
In fact, one has
 \be  L_i = \Lam(X_{ii}) \,. \ee
With a slight abuse of notation, we shall write
 \be  \Lam = (L_1,L_2\ty,\ldots,L_m|L_{m+1},L_{m+2}\ty,\ldots,L_{m+n}) \,.
                                                          \label{labels} \ee
Let us also mention that we are going to use the so--called distinguished 
system of simple roots.

Next we add a few comments on gradations. As is well--known, the algebra
$sl(m|n)$ has a natural $\ZZ$--gradation which is consistent with its
$\ZZ_2$--gradation \cite{S}. Consequently, $U(sl(m|n))$ is $\ZZ$--graded as
well. These $\ZZ$--gradations are easily described by means of the element
 \be  D = - \sum_{j=1}^{m} X_{jj} = \sum_{k=m+1}^{m+n} X_{kk}
  = \frac{1}{m-n}\:\!%
  \mbox{diag}(\underbrace{n,n,\ldots,n}_{m};\underbrace{m,m,\ldots,m}_{n}) \,.
                                                             \label{defD} \ee
In fact, if $X \in sl(m|n)$ or $X \in U(sl(m|n))$, then $X$ is homogeneous of
$\ZZ$--degree $r$ if and only if
 \be  \la D,X \ra = r X \,.  \ee
(Here and in the following $\la\;,\,\ra$ denotes the super commutator.)

On the other hand, if $V$ is a graded $sl(m|n)$--module (in the $\ZZ_2$--graded
or $\ZZ$--graded sense), we can shift the gradation and obtain another graded
$sl(m|n)$--module (see Ref.~\cite{Sgtc} or also Ref.~\cite{SZ}). If $V$ is
finite--dimensional and simple, its gradation can be fixed by specifying the
degree of a highest weight vector. In the following, we shall mainly be
interested in finite--dimensional simple modules $V$ for which the coefficients
$L_i$ of the highest weight $\Lam$ are integral. In this case, there is
a preferred choice of the gradation: It is given by demanding that an element
$x \in V$ be homogeneous of $\ZZ$--degree $r$ if and only if
 \be  D \cd x = r\,x  \label{Dgrad} \ee
(where the dot denotes the action of $D$ on $x$ in $V$\/). Then, for any
finite--dimensional simple graded subquotient $V$ of $U(sl(m|n))$ (endowed
with the adjoint action), the $\ZZ$--gradation of $V$ induced from that of
$U(sl(m|n))$ is exactly the $\ZZ$--gradation specified by Eq.~\reff{Dgrad}.
Thus in the following the gradation will normally be given by that equation.
(The vector module $W$ is an exception: For this module we normally choose
$W_0 = W_{\bar{0}}$ and $W_1 = W_{\bar{1}}$\,, with $\dim W_0 = m$ and
$\dim W_1 = n$\,, however, in Sec.~3 we make a different choice.) Actually,
for our purposes the precise specification of the gradations is, to some
extent, not even necessary: A shift of the gradation of the module of
coefficients simply results in a shift of the gradation of the cohomology
groups \cite{SZ}.

In the following (up to the appendix), we consider the algebras $sl(m|1)$ with
$m \geq 2$ and set
  \be  L = sl(m|1) \,. \label{defL} \ee
We note that in this case we have
  \be D = X_{m+1,\:\!m+1} \,. \ee
(The case of the algebra $sl(2|1) \simeq sl(1|2)$ has been considered in more
detail in Ref.~\cite{SZ}.)

Our goal is to prove Eq.~\reff{claim} for $L = sl(m|1)$. We are going to
proceed as in Ref.~\cite{SZ} and show that 
 \be  H^2(L,V) = \Ze  \label{cVtr} \ee
for all finite--dimensional simple (graded) subquotients $V$ of the $L$--module
$U(L)$.

Let $V = V(\Lam)$ be a finite--dimensional simple $L$--module with
highest weight
 \be  \Lam = (L_1,L_2\ty,\ldots,L_m|L_{m+1}) \,,  \ee
where the $L_i$ satisfy Eq.~\reff{relweight}. According to Ref.~\cite{SZ},
we have
 \be  H^n(L,V) = \Ze \quad\mbox{for all $n$} \ee
whenever $\Lam$ does not belong to one of the following two families of
weights: 
\bea
 \mbox{(0)} & \hspace{3em} & (p,1,\ldots,1\ty|-p-(m-1))
   \hspace{1.5em} \mbox{with $p \geq 1$ integral,\hspace{3em}} \label{famz} \\
 \mbox{(1)} & \hspace{3em} & \hspace{2.5em}(0,\ldots,0,-q\ty|\ty q)
               \hspace{1.5em} \mbox{with $q \geq 0$ integral.} \label{famo}
\eea
(The families are labelled by the number $k \in \{0,1\}$ appearing in their
definition.)
Recall that these are the highest weights of those finite--dimensional simple
$L$--modules for which all Casimir operators without a constant term are
equal to zero. Quite generally, if $V$ is a module of this type, then its dual
is likewise. In the present case, the module of the family (1) with $q = 0$
is trivial and hence self--dual, and it can be shown that the modules of the
families (0) and (1) with $p = q \geq 1$ are dual to each other. In
particular, if $V$ is a simple module of the family (1) with highest weight
$(0,\ldots,0,-p\ty|\ty p)$ (where $p \geq 1$ is an integer), then the
representative $D_V$ of $D$ in $V$ has exactly the eigenvalues
 \be  p,p+1,\ldots,p+m-1 \,, \label{Dev} \ee
and the corresponding eigenvalues in a simple module belonging to the family
(0) are the numbers opposite to those in \reff{Dev}.

On the other hand, the eigenvalues of $D_{U(L)}$ are the numbers
$0,\pm 1,\ldots,\pm m$\,. Thus, from the simple modules in the families (0)
and (1), only the trivial module with highest weight $(0,\ldots,0\ty|\ty 0)$
and the two contragredient modules with the highest weights
$(1,\ldots,1\ty|\ty m)$ resp. $(0,\ldots,0,-1\ty|\ty 1)$ can be isomorphic to
a finite--dimensional simple subquotient of $U(L)$. Consequently, our claim
will be proved if we can show that Eq.~\reff{cVtr} holds for each of these
three modules.

Actually, our task can be simplified even more. As is well--known, the
mapping
 \be  \tau : L \lra L \quad,\quad \tau(A)
                        = -\,{}^{\mbox{\scriptsize{st}}}\!A \label{auto} \ee
(where \mbox{$\,{}^{\mbox{\scriptsize{st}}}\!\!\:A$} is the super transpose
of $A$) is an automorphism of the Lie superalgebra $L$\,. Moreover, if $V$ is a
finite--dimensional simple $L$--module and if $\vr$ is the representation
afforded by $V$, then the representation $\vr \comp \tau$ is equivalent to the
representation contragredient to $\vr$\,. But if $V^{\tau}$ denotes the
graded vector space $V$, endowed with the representation $\vr \comp \tau$\,,
then we know that $H^n(L,V)$ and $H^n(L,V^{\tau})$ are isomorphic (see
Eq.~(2.34) of Ref.~\cite{SZ}). The upshot is that we only have to
consider the trivial module and the simple module with highest weight
$-\ve_m + \ve_{m+1}$\,.
%

\sect{Completion of the proof of our claim}
%
%
According to the previous section we have to show that the 2--cohomology of
$L$ with coefficients in the trivial module $V(0) = \KK$ and in
$V(-\ve_m + \ve_{m+1})$ is trivial.

The case of the trivial $L$--module $\KK$ is easy (see also Ref.~\cite{Fu}).
According to Prop.~2.1 of Ref.~\cite{SZ} every cohomology class in $H^2(L,\KK)$
contains an $L_0$--invariant cocycle. A short look at the representations of
$L_0$ carried by $L_0$ and $L_{\pm 1}$ shows that there exists, up to the
normalization, a unique non--zero super--skew--symmetric $L_0$--invariant
bilinear form on $L$\,. Obviously,
 \be  (A,B) \lra \mbox{Tr}(\la A,B \ra)  \ee
is such a form, but this form is a coboundary. Thus we have shown that
 \be  H^2(L,\KK) = \Ze \,. \ee

The case of the module
 \be  V = V(-\ve_m + \ve_{m+1})  \ee
is more difficult. First of all, we need a suitable realization of this
module. Since we are going to use the Lie superalgebra $gl(m|1)$, we introduce
the abbreviation
 \be  G = gl(m|1) \,. \label{defG} \ee
Let $W$ be the vector module of $G$\,, but endowed with shifted $\ZZ$-- and
$\ZZ_2$--gradations such that
 \be W_{-1} = W_{\bar{1}} \quad,\quad \dim W_{\bar{1}} = m \label{grW1} \ee
 \be  W_{0} = W_{\bar{0}} \quad,\quad \dim W_{\bar{0}} = 1 \label{grW2} \ee
(the reason for this unusual choice will become obvious below, see Remark 3.1),
moreover, let $S(W,\ve)$ be the super--symmetric algebra constructed over $W$
(with $\ve$ the standard commutation factor of supersymmetry; see
Ref.~\cite{Sgtc}). We don't write a product sign for the multiplication in
$S(W,\ve)$. It is well--known that there exists a natural representation
$\vr_0$ of $G$ in $S(W,\ve)$, defined such that, for every $A \in G$\,, the
representative $\vr_0(A)$ is the unique super derivation of $S(W,\ve)$ which
extends $A$\,.

Let $\vr$ denote the representation obtained from $\vr_0$ by a certain twist:
 \be  \vr(A) = \vr_0(A) - \mbox{Str}(A)\ty id  \ee
(where Str denotes the super trace). Obviously, the components $S_n(W,\ve)$
are invariant under $\vr_0$ and $\vr$\,, and it is easy to see that
$S_{m-1}(W,\ve)$, endowed with the representation of $L$ induced by $\vr$\,,
is isomorphic to the $L$--module $V$. This is the realization of $V$ that we
are going to use in the sequel.

More explicitly, let $(\th_i)_{1 \leq i \leq m}$ be a basis of $W_{\bar{1}}$\,,
let $z$ span \vspace{0.2ex}
the one--dimensional space $W_{\bar{0}}$\,, and let
$\frac{\rd}{\rd \th_i}$ and $\frac{\rd}{\rd z}$ be the corresponding
super derivations of $S(W,\ve)$. Then $\vr$ is given by
 \bea
       \vr(E_{i,\ty j}) & = & \th_i {\dis \frac{\rd}{\rd \th_j}} - \de_{ij}
                                                                    \\[0.5ex]
     \vr(E_{i,\ty m+1}) & = & \th_i {\dis \frac{\rd}{\rd z}} \\[0.5ex]
     \vr(E_{m+1,\ty i}) & = & z {\dis \frac{\rd}{\rd \th_i}} \\[0.5ex]
   \vr(E_{m+1,\ty m+1}) & = & z {\dis \frac{\rd}{\rd z}} + 1 \:,
 \eea
where $i,j \in \{1,2\ty,\ldots,m\}$.

\noindent
{\em Remark 3.1.} At this point it should be clear why we have chosen the
gradation of $W$ as given in Eqs.~\reff{grW1},\,\reff{grW2}. With this choice,
the $\th_i$ are fermionic variables, and $z$ is bosonic in the usual sense.
In particular, $z$ commutes with the $\th_i$\,. We could also work with the
standard gradation of $W$. Then $S(W,\ve)$ must be replaced by the
super--Grassmann algebra constructed over $W$, the $\th_i$ are still
fermionic and $z$ is bosonic, however, now $z$ {\em anticommutes} with the
$\th_i$\,. In principle, there is nothing wrong with this choice, but we
wanted to avoid this unusual situation.

Regarded as an $L_0$--module, (in fact, also as a $G_0$--module,)
 \be V = S_{m-1}(W,\ve) \ee
decomposes into
 \be V = \bop^{m}_{r=1} V_r \,,\ee
where
 \be V_r = z^{r-1} \bigwedge^{m-r} W_{\bar{1}} \ee
is a simple $L_0$--module with highest weight
 \be -\ve_{m-r+1} - \ve_{m-r+2} - \ldots - \ve_m + r\:\!\ve_{m+1} \ee
and highest weight vector
 \be z^{r-1} \th_1 \th_2 \ldots \th_{m-r} \,. \ee
Note that according to Eq.~\reff{Dgrad} $V_r$ is the $\ZZ$--homogeneous
component of $V$ of degree $r$\,. The gradation inherited from $W$ is
obtained from this one by a shift.

We now are ready to determine $H^2(L,V)$. Once again by Prop.~2.1 of
Ref.~\cite{SZ}, every cohomology class in $H^2(L,V)$ contains an
$L_0$--invariant cocycle. To find these cocycles, we make a detour via
$G = gl(m|1)$. Let
 \be f : L \ti L \lra V \ee
be any bilinear mapping. Define the bilinear mapping
 \be \bar{f} : G \ti G \lra V \ee
by setting
 \be \bar{f}(A,B) = f(A,B) \hspace{2em}\mbox{for all $A,B \in L$} \ee
 \be \bar{f}(I,C) = -\bar{f}(C,I) = 0 \hspace{2em}\mbox{for all $C\in G$} \ee
(recall that $I$ denotes the $(m+1) \ti (m+1)$ unit matrix). Then it is easy
to see that $f$ is an $L_0$--invariant 2--cocycle if and only if $\bar{f}$ is
a $G_0$--invariant 2--cocycle (note that $I$ acts on $V$ by the zero operator).
Consequently, it is sufficient to determine the $G_0$--invariant 2--cocycles
on $G$ (with values in $V$\/); the $L_0$--invariant 2--cocycles on $L$
are then simply obtained by restriction.

Let
 \be g : G \ti G \lra V \ee
be a super--skew--symmetric $G_0$--invariant bilinear mapping. According to
our conventions, the invariance of $g$ under $D$ says that $g$ is homogeneous
of degree $0$ in the sense of $\ZZ$--gradations. Thus if $r$ and $s$ are two
elements of $\{-1,0,1\}$, it follows that
 \be g(G_r \ti G_s) \subset V_{r+s} \ee
and hence that
 \be g(G_r \ti G_s) = \Ze \quad\mbox{if $r+s \leq 0$} \,. \ee
Moreover, the restriction of $g$ to $G_1 \ti G_1$ must be symmetric. But
a look at the $sl(m)$--module structures of $S_2(G_1)$ (the symmetric tensor
product of $G_1$ with itself) and $V_2$ shows that a non--zero symmetric
$sl(m)$--invariant bilinear mapping of $G_1 \ti G_1$ into $V_2$ does not
exist. Thus we have
 \be g(G_1 \ti G_1) = \Ze \,, \ee
and all we have to do is to find the restriction of $g$ onto $G_0 \ti G_1$\,,
say.

To construct the $G_0$--invariant bilinear mappings $G_0 \ti G_1 \rar V_1$\,,
we first define, for $i \in \{1,2\ty,\ldots,m\}$,
 \be \eta_i = (-1)^{i-1} \th_1 \ldots {\widehat{\th}}_i \ldots \th_m
                      = \frac{\rd}{\rd \th_i} (\th_1 \th_2 \ldots \th_m) \ee
(as usual, the sign \,$\what{0.5ex}$\, indicates that the element below it must
be omitted). Obviously, the $\eta_i$ form a basis of $V_1$\,, moreover, we have
 \be \vr(E_{ij})\ty \eta_k = -\de_{ik}\ty \eta_j \ee
for all $i,j,k \in \{1,2\ty,\ldots,m\}$. Combined with the known actions of
$D$ and $I$\,, this shows explicitly that the $G_0$--modules $G_1$ and $V_1$
are isomorphic.

Using this information as well as the standard representation theory of
$sl(m)$, we now can describe the super--skew--symmetric $G_0$--invariant
bilinear mappings \\
$G \ti G \rar V$, as follows. Define three bilinear maps
$g_1,g_2\ty,g_3$ of $G \ti G$ into $V$ by
 \be g_1(E_{i,\ty j}\ty ,E_{m+1,\ty k})
               = -g_1(E_{m+1,\ty k}\ty ,E_{i,\ty j}) = \de_{ik}\ty \eta_j \ee
 \be g_2(E_{i,\ty j}\ty ,E_{m+1,\ty k})
               = -g_2(E_{m+1,\ty k}\ty ,E_{i,\ty j}) = \de_{ij}\ty \eta_k
                                                           \vspace{0.8ex} \ee
 \be g_3(E_{m+1,\ty m+1},E_{m+1,\ty k})
               = -g_3(E_{m+1,\ty k}\ty ,E_{m+1,\ty m+1}) = \eta_k \,, \ee
where $i,j,k \in \{1,2\ty,\ldots,m\}$, and with the understanding that the
values of $g_1,g_2\ty,g_3$ on the remaining pairs of the standard basis
elements of $G$ are equal to zero. Then $g_1,g_2\ty,g_3$ are
super--skew--symmetric and $G_0$--invariant, and any bilinear map
$g : G \ti G \rar V$ with these properties is a linear combination of them.

Now suppose that $g : G \ti G \rar V$ is a $G_0$--invariant 2--cocycle. Using
the $G_0$--invariance of $g$ as well as the fact that $g$ vanishes on
$G_0 \ti G_0$\,, it is easy to see that
 \be (\de^2 g)(A,B,C) = g(\la A,B \ra, C) \ee
for all $A,B \in G_0$ and $C \in G$\,. Hence the cocycle condition
implies that $g$ vanishes on $sl(m) \ti G$\,. Consequently, $g$ must be a
linear combination of $g_2$ and $g_3$\,,
 \be g = a\ty g_2 + b\ty g_3 \,. \ee
A short calculation then shows that
 \be (\de^2 g)(E_{k,\ty m+1},E_{m+1,\ty i}\ty,E_{m+1,\ty j})
                     = -(a + b)(\de_{ki}\ty \eta_j  + \de_{kj}\ty \eta_i) \ee
for all $i,j,k \in \{1,2\ty,\ldots,m\}$, which implies that
 \be a + b = 0 \,. \ee
Without loss of generality we now may assume that
 \be a = 1 \,, \ee
and then we have
 \be g(A,E_{m+1,\ty k}) = \mbox{Str}(A)\ty \eta_k \ee
for all $A \in G$ and all $k \in \{1,2\ty,\ldots,m\}$. Consequently, $g$
vanishes on $L \ti L$\,, and according to our previous discussion, this
implies that
 \be H^2(L,V) = \Ze \,, \ee
as claimed.

We close this section by the remark that $g$ as specified above is a 2--cocycle
on $G$\,, and that $g$ is not a 2--coboundary. Thus we have
 \be \dim H^2(G,V) = 1 \,. \ee
%

\sect{Discussion}
%
%
In the present paper we have shown that
 \be H^2(L,U(L)) = \Ze \label{cl} \ee
for the Lie superalgebras $L = sl(m|1)$ with $m \geq 2$\,. Our method of
proof was the following.

Because of the long exact cohomology sequence \cite{SZ}, a sufficient (but not
necessary) condition for \reff{cl} to hold is that
 \be H^2(L,V) = \Ze \label{clv} \ee
for all simple subquotients $V$ of $U(L)$ (these are automatically
finite--dimensional).

Let $\Lam$ be the highest weight of $V$. To prove that Eq.~\reff{clv} holds
for the modules $V$ in question, we first used the results of Example 4 in
Sec.~3 of Ref.~\cite{SZ} (and hence Prop.~2.2 of that reference) to conclude
that Eq.~\reff{clv} is true if $\Lam$ does not belong to the families (0) and
(1) defined by Eqs.~\reff{famz},\,\reff{famo}. By comparing the eigenvalues of
$D_{U(L)}$ and $D_V$ we could then reduce the problem to the consideration of
just three cases, finally, by using the automorphism \reff{auto} of $L$\,,
even to two cases. One of these cases was the trivial module, for which
Eq.~\reff{clv} could be proved immediately. In the other case, we found a
nice realization of $V$, which made the necessary calculations simple.

In view of our experience with $sl(3|2)$ (see the appendix) it must be said
that the case of the algebras $sl(m|1)$ is particularly favourable. In more
general cases (already for $sl(3|2)$\ty) we certainly need more information on
the adjoint $L$--module $U(L)$ than just the eigenvalues of $D_{U(L)}$\,. Also,
most of the modules $V$ one finally has to consider will not be well--known,
and it may be very hard to find a suitable realization for them. (Recall that,
at least for the $sl(m|n)$ algebras with $m \neq n$\,, these modules are
maximally atypical \cite{SZ}.)

All this seems to indicate that in our approach too many details are needed,
and that more profound methods are necessary to solve our problem for more
general algebras.

\vspace{2ex}

\noindent
{\bf Acknowledgement} \\
The present work was initiated during a visit of the first--named author to
the Department of Pure Mathematics of the University of Adelaide. The kind
invitation by the second--named author and the hospitality extended to the
first--named author, both in the Mathematics and in the Physics Department,
are gratefully acknowledged.

\vspace{8ex}

\noindent
{\LARGE\bf Appendix}\\[-6ex]

\begin{appendix}
%

\sect{Problems with $sl(3|2)$}
%
%
Once we had proved Eq.~\reff{claim} for the algebras $L = sl(m|1)$\,;
$m\geq 2$, we intended to investigate the algebras $sl(m|n)$\,; $m \neq n$\,,
in general. In order to see what type of problems would arise, we first
considered the algebra $sl(3|2)$. Unfortunately, up to now we have not been
able to prove or disprove Eq.~\reff{claim} in this case. Nevertheless, we
think it may be worthwhile to present some of our intermediate results, in
order to show the obstacles one has to overcome if one wants to proceed
along the lines described in Ref.~\cite{SZ} and in the present paper.

In this appendix, we set
 \be L = sl(3|2) \ee
and use the notation introduced at the beginning of Sec.~2. We try to proceed
as in Secs.~2 and 3. Let $V(\Lam)$ be a finite--dimensional simple $L$--module
with highest weight
 \be \Lam = (L_1,L_2\ty,L_3|L_4,L_5) \,, \ee
where the $L_i$ satisfy Eq.~\reff{relweight}. According to Ref.~\cite{SZ}, we
have
 \be H^n(L,V(\Lam)) = \Ze \quad\mbox{for all $n$} \ee
whenever $\Lam$ does not belong to one of the following three families of
weights:
\bea
  \mbox{(0)} & \hspace{3em} & (p,q,2\ty|-q-1,-p-1)
               \hspace{1.5em} \mbox{with $p \geq q \geq 2$\hspace{3em}}
                                                            \label{faone} \\
  \mbox{(1)} & \hspace{3em} & \hspace{1.0em}(p,1,q\ty|-q,-p-1)
               \hspace{1.5em} \mbox{with $p \geq 1 \geq q$}
                                                            \label{fatwo} \\
  \mbox{(2)} & \hspace{3em} & \hspace{2.0em}(0,p,q\ty|-q,-p)
               \hspace{1.5em} \mbox{with $0 \geq p \geq q$}
                                                          \label{fathree}
\eea
where, in all cases, $p$ and $q$ are {\em integers.}

We already know that the class of simple $L$--modules with these highest
weights is closed under taking duals. The following list gives, on the same
line, for each of the weights appearing in Eqs.~\reff{faone} -- \reff{fathree},
the highest weight of the corresponding dual module: \addtocounter{equation}{1}
\[\begin{array}{@{}c@{\hspace*{1.0em}}r@{\hspace{0.6em},\hspace{0.6em}}%
l@{\hspace{1.0em}}l@{\hspace{1.14em}}r@{}} 
   & (p,q,2\ty|-q-1,-p-1) & (0,1-q,1-p\ty|\ty p-1,q-1)
        & \mbox{with $p > q \geq 2$} & \mbox{(A.\arabic{equation})}
                                       \addtocounter{equation}{1} \\[0.5ex]
   & (p,p,2\ty|-p-1,-p-1) & (0,-p,-p\ty|\ty p,p)
        & \mbox{with $p \geq 2$}     & \mbox{(A.\arabic{equation})}
                                       \addtocounter{equation}{1} \\[0.5ex]
   & (p,1,q\ty|-q,-p-1)   & (1-q,1,1-p\ty|\ty p-1,q-2)
        & \mbox{with $p \geq 1 > q$} & \mbox{(A.\arabic{equation})}
                                       \addtocounter{equation}{1} \\[0.5ex]
   & (p,1,1\ty|-1,-p-1)   & (0,0,1-p\ty|\ty p-1,0)
        & \mbox{with $p \geq 2$}     & \mbox{(A.\arabic{equation})}
                                       \addtocounter{equation}{1} \\[0.5ex]
   & (1,1,1\ty|\ty -1,-2) & (0,-1,-1\ty|\ty 1,1)
        &                            & \mbox{(A.\arabic{equation})}
                                       \addtocounter{equation}{1} \\[0.5ex]
   & (0,0,0\ty|\ty 0,0)   & (0,0,0\ty|\ty 0,0) \;.
        &                            & \mbox{(A.\arabic{equation})}
  \end{array} \]
We note that, in this list, a weight appears twice if and only if the
corresponding module is self--dual. Let us also mention that, apart from the
trivial module, every module of the class (2) is dual to some of the modules
in the classes (0) or (1). Similarly, every module of the class (0) is
dual to some module of the class (2).

Next we note that the representative of $D$ under the adjoint representation
in $U(L)$ has exactly the eigenvalues $0,\pm 1,\ldots,\pm 6$\,. Consequently,
for any subquotient $V$ of $U(L)$, the eigenvalues of the representative
$D_V$ must belong to this set. To apply this condition, we note that, for any
weight $\Lam$ given by the Eqs.~\reff{defweight},\,\reff{relweight}, we have
 \be \Lam(D) = -L_1 - L_2 - L_3 = L_4 + L_5 \,. \ee

At this point, we meet a first complication. Whereas only a few modules of the
types (0) and (2) pass the criterion above, there are infinitely many modules
of type (1) which satisfy it. Thus we have to find some other conditions to
narrow down the number of possibilities. Unfortunately, we don't have any more
detailed information on the structure of the $L$--module $U(L)$. Consequently,
for the time being, we simply ignore the fact that the modules we are finally
interested in are simple subquotients of $U(L)$. Instead, we consider all the
simple $L$--modules $V(\Lam)$, with $\Lam$ a weight of the types
(0),\,(1),\,(2) above, and we try to find a simple necessary criterion
implying that $H^2(L,V(\Lam)) = \Ze$\,.

This can be achieved as follows. Once again because of Prop.~2.1 of
Ref.~\cite{SZ}, we are only interested in the $L_0$--invariant 2--cocycles
with values in $V(\Lam)$. In particular, such a cocycle can be identified with
an $L_0$--module homomorphism of $L \wepf L$ (the super--exterior square of
the adjoint $L$--module) into $V(\Lam)$. But a non--zero homomorphism of this
type exists if and only if there is at least one simple $L_0$--submodule of
$L \wepf L$ which is isomorphic to an $L_0$--submodule of $V(\Lam)$.

The $L_0$--module structure of $L \wepf L$ is easily determined: Suffice it
to say that $L \wepf L$ decomposes into the direct sum of 27 simple
$L_0$--submodules (not all non--isomorphic, of course).

The $L_0$--module structure of the modules $V(\Lam)$ is much more difficult
to obtain. Information on these modules could be extracted from the
conjectured character formula in Ref.~\cite{VHKT} or from the character formula
proved in Ref.~\cite{Ser}. However, we haven't tried to do that but rather
argue more directly, as follows.

Let $\bar{V}(\Lam)$ be the Kac module with highest weight $\Lam$ (see
Ref.~\cite{Kabo}). If there is a simple $L_0$--submodule of $L \wepf L$
which is isomorphic to a submodule of $V(\Lam)$, this is even more the case
with $V(\Lam)$ replaced by $\bar{V}(\Lam)$ (since $V(\Lam)$ is a quotient of
$\bar{V}(\Lam)$\ty). But the structure of the $L_0$--module $\bar{V}(\Lam)$
can be determined by the standard representation theory of $L_0$\,, and then
it is not difficult to check whether the condition above is satisfied. (All
this is a straightforward but cumbersome task: Note that $\bar{V}(\Lam)$ may
be the direct sum of up to 64 simple $L_0$--submodules.)

The calculation sketched above yields a finite list of highest weights
$\Lam$\,. But there is still one more observation to be made. Obviously, the
$L$--module $L \wepf L$ is self--dual. Thus if there exists a non--trivial
$L_0$--module homomorphism of $L \wepf L$ into $V(\Lam)$, then there is also
one into the dual module $V(\Lam)^{\ast}$. Since we have been working with
the Kac modules $\bar{V}(\Lam)$, there may be --- and indeed are -- weights
$\Lam$ in the list above for which the highest weight of the dual module
$V(\Lam)^{\ast}$ is not in the list. All these weights $\Lam$ may be excluded
as well.

The upshot of all this is that, for any finite--dimensional simple $L$--module
$V(\Lam)$ with highest weight $\Lam$\,, the inequality
$H^2(L,V(\Lam)) \neq \Ze$ implies that $\Lam$ is one of the following
weights:
 \be (0,0,0\ty|\ty 0,0) \label{wtriv} \ee
 \be (1,1,1\ty|-1,-2) \;\;,\;\; (0,-1,-1\ty|\ty 1,1) \label{wvec} \ee
 \be (2,1,1\ty|-1,-3) \;\;,\;\; (0,0,-1\ty|\ty 1,0)  \ee
 \be (1,1,0\ty|\ty 0,-2) \label{wstwo} \ee
 \be (3,1,1\ty|-1,-4) \;\;,\;\; (0,0,-2\ty|\ty 2,0)  \ee
 \be (2,1,0\ty|\ty 0,-3) \;\;,\;\; (1,1,-1\ty|\ty 1,-2) \label{wlast} \ee
and, in addition,
 \be (3,1,0\ty|\ty 0,-4) \;\;,\;\; (1,1,-2\ty|\ty 2,-2) \label{addone} \ee
 \be (2,1,-1\ty|\ty 1,-3) \,. \label{addtwo} \ee
As before, if two weights stand on the same line, the corresponding
$L$--modules are dual to each other; if there is only one weight on a line, the
corresponding $L$--module is self--dual. The reader can easily convince
himself/herself that, for each of these weights, the representative
$D_{V(\Lam)}$ takes its eigenvalues in the set $\{-6,-5,\ldots,6\}$. Hence
the assumption that $V(\Lam)$ is isomorphic to a subquotient of $U(L)$
doesn't imply any further restrictions.

Our next task would be to calculate $H^2(L,V(\Lam))$ for the weights $\Lam$
given above. To do this we need much more information on the modules
$V(\Lam)$. In particular, we need the $L_0$--module structure of these
modules. In a painstaking analysis of the corresponding Kac modules
$\bar{V}(\Lam)$, we have determined how the modules $V(\Lam)$ decompose when
regarded as $L_0$--modules. As a by--product, we have also found the
composition factors of the Kac modules themselves. It turns out that the
$V(\Lam)$ in question admit a unique decomposition into simple
$L_0$--submodules, i.e., all the simple $L_0$--modules contained in $V(\Lam)$
have multiplicity one.

The more detailed information thus obtained allows us to rule out the weights
in Eqs.~\reff{addone},\,\reff{addtwo}: For these weights, $V(\Lam)$ doesn't
contain a simple $L_0$--submodule which is isomorphic to an $L_0$--submodule 
of $L \wepf L$\,. Thus we are left with the weights in
Eqs.~\reff{wtriv} -- \reff{wlast}. Using the automorphism \reff{auto} for
$L = sl(3|2)$, we only have to consider one of the weights on each line.
Thus, there are six cases to consider.

The first case is the trivial $L$--module $\KK$ with the highest weight
\reff{wtriv}. It is known from Ref.~\cite{Fu} and easy to see by means of
Prop.~2.1 in Ref.~\cite{SZ} that
 \be H^2(L,\KK) = \Ze \,. \ee
The weights in \reff{wvec} correspond to the covector and vector modules of
$L$\,. Somewhat unexpectedly, for these modules $V$ we have
 \be \dim H^2(L,V) = 1 \,. \ee
This is bad news: It shows that in order to prove Eq.~\reff{claim} for
$L = sl(3|2)$ we need more detailed information on the structure of the
$L$--module $U(L)$.

At this point, we have changed our strategy: Maybe $H^2(L,U(L))$ is different
from $\Ze$. In order to show this we recall that $U(L)$, regarded as an
$L$--module, is canonically isomorphic to the super--symmetric algebra
$S(L,\ve)$ (see Ref.~\cite{S.ce}). It is well--known that $S(L,\ve)$ decomposes
into the direct sum of its $\ZZ$--homogeneous components $S_n(L,\ve)$\,;
$n \geq 0$\,, and that these are $L$--submodules of $S(L,\ve)$, moreover,
$S_n(L,\ve)$ is canonically isomorphic to the submodule of super--symmetric
tensors in $L^{\ot n}$. In particular, the submodule $S_0(L,\ve)$ is isomorphic
to $\KK$\,, and we already know that $H^2(L,\KK)$ is trivial. The submodule
$S_1(L,\ve)$ is isomorphic to the adjoint $L$--module, and its highest weight
$(1,0,0\ty|\ty 0,-1)$ is not of type (0), (1), or (2). Hence $H^2(L,L)$ is
trivial as well.

Thus we consider $S_2(L,\ve)$. A detailed analysis shows that the $L$--module
$S_2(L,\ve)$ (uniquely) decomposes like
 \be S_2(L,\ve) \simeq V(2,0,0\ty|-1,-1) \op V(1,0,0\ty|\ty 0,-1) \op W \,, \ee
with an indecomposable but non--simple $L$--module $W$. The module $W$ has a
Jordan--H\"older series
 \be W = W_0 \supset W_1 \supset W_2 \supset W_3 = \Ze \ee
such that $W/W_1$ and $W_2$ are trivial one--dimensional $L$--modules
and such that
 \be W_1/W_2 \simeq V(1,1,0\ty|\ty 0,-2) \,. \ee
Note that the latter module has the highest weight given in \reff{wstwo}.

The module $W_2$ consists of the $L$--invariant elements in $S_2(L,\ve)$,
it is spanned by the so--called (quadratic) tensor Casimir (split Casimir)
element $C$\,. (It is known that any $L$--invariant element of $S_2(L,\ve)$ is
proportional to $C$\,: Otherwise, $L$ would have two linearly independent
quadratic Casimir elements or, equivalently, two linearly independent
super--symmetric $L$--invariant bilinear forms.)

On the other hand, let $G$ be an $L_0$--invariant element of $W$ which does not
belong to $W_1$\,. Any other element $G'$ with these properties has the form
 \be G' = a\ty G + b\ty\ty C \ee
with a non--zero constant $a$ and an arbitrary constant $b$\,. According to
the preceding characterization of $W_2$\,, the element $G$ is {\em not}
$L$--invariant. Actually, $G$ generates the $L$--module $W$ (but, of course,
it is not a highest weight vector). Let us also note that $W_1$ is the unique
maximal and $W_2$ the unique minimal (i.e., simple) $L$--submodule of $W$.

\noindent
{\em Remark A.1.} Obviously, the module $S_2(L,\ve)$ is self--dual, and so
are the modules $V(2,0,0\ty|-1,-1)$ and $V(1,0,0\ty|\ty 0,-1)$ (the latter is
isomorphic to the adjoint $L$--module). A moment's thought then shows that $W$
is self--dual as well. This explains part of the structure of $W$. Note that
a similar but even more complicated structure also exists for $sl(1|2)$ (see
Eq.~(B.3) of Ref.~\cite{SZ}).

The weights $(2,0,0\ty|-1,-1)$ and $(1,0,0\ty|\ty 0,-1)$ do not belong
to the families (0), (1),\,(2). Thus we have
 \be H^n(L,S_2(L,\ve)) \simeq H^n(L,W) \quad\mbox{for all $n$} \,, \ee
and a rather tedious calculation shows that
 \be H^n(L,W) = \Ze \quad\mbox{for $n = 1,2$} \,. \ee
As is easily guessed, a lot of detailed knowledge about the modules
$V(1,1,0\ty|\ty 0,-2)$ and $W$ is necessary to prove this result.

Thus the situation is rather unpleasant: We know of simple $L$--modules (the
vector and covector modules) for which the 2--cohomology is non--trivial, but
since we don't have sufficient information on the $L$--module $U(L)$, we do
not know what this fact implies for $H^2(L,U(L))$. On the other hand, one of
the candidates for a non--trivial cohomology (namely $V(1,1,0\ty|\ty 0,-2)$\ty)
really {\em is} isomorphic to a subquotient of $U(L)$, but this does not imply
that $H^2(L,U(L))$ is non--trivial. (We stress that we have {\em not} shown
that the 2--cohomology of $L$ with values in $V(1,1,0\ty|\ty 0,-2)$ is
trivial.)

\end{appendix}

\end{document}